\documentclass[12pt, a4paper]{article}
\usepackage[cp1251]{inputenc}
\usepackage{graphics}
\usepackage{amsmath}
\usepackage{amssymb}
\usepackage{enumitem}
\usepackage{amsfonts}
\usepackage{geometry}
\usepackage{amsthm, latexsym}
\usepackage[english]{babel}
\newtheorem{theorem}{\qquad Theorem}
\newtheorem{proposition}{\qquad Proposition}

\newtheorem{lemma}{\qquad Lemma}

\title{On the Riemann-Hilbert Problem for Difference and $q$-Difference Systems}

\author{Ilya Vyugin, \; Roman Levin}

\date{}

\begin{document}

\maketitle


\bigskip


\begin{abstract}
In this paper we study an analogue of the classical Riemann-Hilbert problem stated for the classes of difference and $q$-difference systems. The Birkhoff's existence theorem was generalized in this paper.
\end{abstract}
\section{Introduction}
The analytic theory of difference and $q$-difference equations was introduced at the beginning of the 20-th century and was completely developed by George D. Birkhoff. G.D.Birkhoff introduced formal and genuine solutions to the difference and $q$-difference linear systems, he also defined a specific periodic ($q$-periodic) matrix function which he considered as a concept of the monodromy for those systems. We have to notice that that matrix function -- the monodromy of difference systems, being defined not in the sense of the solutions branching, but as a ratio of two special fundamental matrices of the systems, bears a remarkable resemblance to the Stokes matrices in the theory of linear differential equations.\\
\\
Having found exact forms of the monodromy and the coefficient matrix of the system, Birkhoff formulated an analogue of the Riemann-Hilbert problem for the difference and $q$-difference cases which he called a generalized Riemann-Hilbert problem. To prevent a confusion, we have to emphasize that there exist several more problems, different from the one we consider in this paper, but also referred to as generalized Riemann-Hilbert problems. The problem which was stated by Birkhoff and to which he proposed a solution could be formulated as follows:\\

\textit{Construct a difference ($q$-difference) linear system with a given monodromy matrix, a prescribed set of characteristic constants, and a condition for the coefficient matrix to be a polynomial of the fixed power.} \\
\\
However, the result obtained by Birkhoff could be insufficient in some cases since his theorem sometimes leads to systems with shifted characteristic constants. There might be integer additions to the characteristic constants corresponding to power asymptotics of the solutions of the systems. In this paper we consider further research on this problem and propose a solution which shows that there exist systems with the correct monodromy data and characteristic constants. As a result, the roots of the determinant of coefficient matrix could be shifted by an integer, but those are not fixed in the monodromy data. 

\section{Difference Systems Case}
\subsection{Introduction to linear difference systems}

A system of linear difference equations is a system of the form:
\begin{eqnarray}\label{dif-syst}
Y(z+1)=A(z)Y(z),
\end{eqnarray}
here $Y(z)$ is a $n\times n$ matrix and $A(z)$ is a rational $n\times n$ matrix of coefficients.  
It could be transformed to the polynomial form 
\begin{eqnarray}\label{syst-koeff}
A(z)=A_rz^r+\ldots+ A_0
\end{eqnarray}
by means of the following gauge transformation:
$$
\tilde{Y}(z)=\Gamma(z-x_1)\cdot\ldots\cdot\Gamma(z-x_s)Y(z),
$$
where $\Gamma(z)$ is the gamma function, $(z-x_1)\cdot\ldots\cdot(z-x_s)$ is the common denominator of elements of the matrix $A(z)$. 
\\
\\
Let us suppose that $\rho_1,\ldots,\rho_n$ are eigenvalues of the matrix $A_r$ and  $\rho_1\cdot\ldots\cdot\rho_n\not=0$, $\rho_i/\rho_j\notin\mathbb{R}$ for $i\not= j$. Then, according to existence theorems (see \cite{Bir-gen}) for difference equations, the formal solution of the system (\ref{dif-syst}) is the following: 
\begin{eqnarray}\label{formal-sol}
\hat{Y}(z)=z^{rz}e^{-rz}\left(\hat{Y}_0+\frac{\hat{Y}_1}{z}+\ldots \right){\rm diag}(\rho_1^z z^{d_1},\ldots,\rho_n^{z}z^{d_n}).
\end{eqnarray}
Without loss of generality, we can suppose that $A_r={\rm diag}(\rho_1,\ldots,\rho_n)$.
It is easy to prove by substitution of $\hat{Y}(z)$ to the system (\ref{dif-syst}); moreover, we obtain that $\rho_1 d_1,\ldots,\rho_n d_n$ coincide with diagonal elements of the matrix $A_{r-1}$ of (\ref{syst-koeff}).\\
\\
Let us denote the roots of the polynomial $\det A(z)$ by $q_1,\ldots,q_{rn}$. It is easy to see that the following holds:
\begin{eqnarray*}
\sum_{i=1}^n d_i+\sum_{j=1}^{rn}q_j=0.
\end{eqnarray*}
This is an analogue of Fuchs relation for linear differential equations.\\
\\

The series in (\ref{formal-sol}) is formal and it could be verified by direct substitution to (\ref{dif-syst}) that it is the solution. The genuine solutions were introduced by Birkhoff. \\

\begin{theorem}[Birkhoff \cite{Bir-gen}, Th. III ]\label{th-formal}
Let us suppose that
$$
A_r={\rm diag}(\rho_1,\ldots,\rho_n),
$$
$$
\rho_1\cdot\ldots\cdot\rho_n\not=0;\qquad\forall i\not=j:\quad \rho_i/\rho_j\notin\mathbb{R}.
$$ 
Then there exists a unique solution $Y^l(z)$ ($Y^r(z)$) of the system (\ref{dif-syst}), such that:
\begin{enumerate}[label = (\roman*)]
\item The function  $Y^l(z)$ ($Y^r(z)$) is analytic in $\mathbb{C}\setminus\{ p_i \}$, where $p_i$ are points on the left (on the right) from poles of $A(z)$ ($A^{-1}(z-1)$) and congruent to them (two points are {\it congruent} if the difference between them is integer);

\item In an arbitrary left (right) half-plain the solution $Y^l(z)$ ($Y^r(z)$) has the asymptotic expansion (\ref{formal-sol}).\\
\end{enumerate}
\end{theorem}
We say that $Y^l(z)$ ($Y^r(z)$) has the {\it asymptotic expansion} (\ref{formal-sol}) in an arbitrary left (right) half-plain if
$$
\left| Y^{l,r}(z)z^{-rz}e^{rz}{\rm diag}(\rho_1^{-z}z^{-d_1},\ldots,\rho_n^{-z}z^{-d_n})-\hat{Y}_0-\frac{\hat{Y}_1}{z}-\ldots-\frac{\hat{Y}_{k-1}}{z^{k-1}}\right|\leqslant \frac{{\rm const}}{z^k},\quad {\rm Re}\, z\to \mp\infty
$$
and the imaginary part of $z$ is bounded.
\subsection{Monodromy matrix of the difference system}

Let us consider the matrix $P(z)=(Y^r(z))^{-1}Y^l(z)$. This matrix $P(z)$ is periodic: 
$$
P(z+1)=P(z).
$$ 
We call the matrix $P(z)$ {\it the monodromy matrix} of the system (\ref{dif-syst}). The exact form of the elements of $P(z)$ is given by the following theorem.\\
\begin{theorem}[Birkhoff \cite{Bir-gen}, Th. IV]\label{th-monodr}
In the assumptions of Theorem \ref{th-formal} the elements $p_{kl}(z)$ of the matrix $P(z)=(Y^r(z))^{-1}Y^l(z)$ can be represented as follows:
\begin{eqnarray}\label{mon-differ}
p_{kk}(z)=1+c_{kk}^{(1)}e^{2\pi i z}+\ldots+c_{kk}^{(r-1)}e^{2\pi (r-1)iz}+e^{2\pi i d_k}e^{2\pi niz},
\end{eqnarray}
$$
p_{kl}(z)=e^{2\pi \lambda_{kl} z}\left( c_{kl}^{(0)}+c_{kl}^{(1)}e^{2\pi i z}+\ldots+c_{kl}^{(r-1)}e^{2\pi (r-1)iz}\right),\quad k\not= l,
$$
where $c_{kl}^{(s)}$ are some constants, and $\lambda_{kl}$ is the minimum integer exceeding ${\rm Re}\frac{\ln \rho_l-\ln\rho_k}{2\pi i}$ for each $k$ and $l$ (we fix branches of $\ln z$ in the left and the right half-plains).\\
\end{theorem}

\subsection{Generalized Riemann-Hilbert problem}
Actually, the matrix polynomial $A(z)=A_rz^r+\ldots+A_0$ with $A_r={\rm diag}(\rho_1,\ldots,\rho_n)$, $\rho_1\cdot\ldots\cdot\rho_n\not=0$, $\rho_k\not=\rho_l$, $k\not=l$ gives us (by Theorems \ref{th-formal} and \ref{th-monodr}) the {\it characteristic constants} $\{ d_k\},\{ c_{kl}^{(s)}\}$.
The number of characteristic constants is equal to the number of elements of the matrices $A_0,A_1,\ldots,A_r$. We will study the map
\begin{eqnarray}\label{Poincare-map}
(A_1,\ldots,A_r)\longmapsto \left(\{ d_k\},\{ c_{kl}^{(s)}\}\right),
\end{eqnarray}
if the constants $\rho_1,\ldots,\rho_n$ are fixed. We can now formulate the generalized Riemann-Hilbert problem for difference systems:\\

{\it Construct a system (\ref{dif-syst}), (\ref{syst-koeff}) with a prescribed set of characteristic constants $\{ d_k\},\{ c_{kl}^{(s)}\}$ and with a given matrix $A_r$ (a reversibility of the map (\ref{Poincare-map}) is studied).}\\ 
\\
Birkhoff formulated the following results:\\

\begin{theorem}[Birkhoff \cite{Bir-Rim}]\label{th-inversprob}
For any nonzero $\rho_1,\ldots,\rho_n$ such that $\forall i\not=j: \rho_i/\rho_j\notin\mathbb{R}$ there exist matrices $A_0,\ldots, A_{r-1}$, such that the system (\ref{dif-syst}), (\ref{syst-koeff}) has the coefficient matrix $A(z)$ with $A_r={\rm diag}(\rho_1,\ldots,\rho_n)$ and has the given characteristic constants $\{ d_k\},\{ c_{kl}^{(s)}\}$ or constants $\{ d_k+l_k\},\{ c_{kl}^{(s)}\}$, where $l_1,\ldots,l_n\in\mathbb{Z}$.\\
\end{theorem}

\begin{theorem}[Birkhoff \cite{Bir-gen}, Th. VII]\label{th-ration-rel}
Let us suppose that there are two matrix polynomials $A'(z)=A_r'z^r+\ldots+A_0'$ and $A''(z)=A_r''z^r+\ldots+A_0''$ with 
$$
A_r'=A_r''={\rm diag}(\rho_1,\ldots,\rho_n),\qquad \rho_1\cdot\ldots\cdot\rho_n\not=0;\quad \rho_k/\rho_l\notin\mathbb{R},\quad k\not= l,
$$
such that the sets of characteristic constants for systems
$$
Y'(z+1)=A'(z)Y'(z),\qquad Y''(z+1)=A''(z)Y''(z)
$$ 
coincide. Then there exists the rational matrix $R(z)$, such that
\begin{eqnarray}\label{ratio-trans}
A''(z)=R(z+1)A'(z)R^{-1}(z),
\end{eqnarray}
and $(Y'')^{l,r}=R(Y')^{l,r}$.\\
\end{theorem}

We propose an improvement to Theorem \ref{th-inversprob} of Birkhoff where we show that it is possible to avoid the shifts of characteristic constants by integers. We now formulate the theorem, the proof will be given later, after the consideration of several preliminary lemmas.

\begin{theorem}\label{th-BSF}
For any nonzero $\rho_1,\ldots,\rho_n$ such that $\rho_i/\rho_j\notin\mathbb{R}$ for $i\not= j$, $\prod_{i=1}^n\rho_i\not=0$ and characteristic constants $\{ d_k\},\{ c_{kl}^{(s)}\}$ such that the matrix $P(z)$ (\ref{mon-differ}) does not have multiple zeros there exist matrices $A_0,\ldots, A_{r-1}$, such that the system (\ref{dif-syst}), (\ref{syst-koeff}) has the coefficient matrix $A(z)$ with  $A_r={\rm diag}(\rho_1,\ldots,\rho_n)$. \\
\end{theorem}

\subsection{Lemmas and proofs}

\begin{lemma}\label{lemma-fixval}
Let us consider 
some system (\ref{dif-syst}), (\ref{syst-koeff}) with coefficients
\begin{eqnarray}\label{cond-syst-1}
A_r={\rm diag}(\rho_1,\ldots,\rho_n),\qquad
\rho_1\cdot\ldots\cdot\rho_n\not=0;\qquad\forall i\not=j:\quad \rho_i/\rho_j\notin\mathbb{R},
\end{eqnarray}
and a set of constants $d_1,\ldots,d_n$. Then for an arbitrary set of parameters $\tilde{d}_i$, $i=1,\ldots,n$, such that $\tilde{d}_i-d_i\in\mathbb{Z}$ for $i=1,\ldots,n$, there exists another system
\begin{eqnarray}\label{syst-tilde}
Y'(z+1) = A'(z) Y'(z)
\end{eqnarray}
with coefficient matrix
\begin{eqnarray}\label{syst-koeff-tilde}
A'(z)=A_r' z^r+A_{r-1}'z^{r-1}+\ldots+A_0'+\ldots+A_{-s}'z^{-s},\quad A_r'=A_r,
\end{eqnarray}
and the characteristic constants $\tilde{d}_i$, $i=1,\ldots,n$, such that $ 
A'(z)=M(z+1)A(z)M^{-1}(z)$, where $M(z)$ is a rational matrix.\\
\end{lemma}

{\it Proof.} Let us prove the lemma by induction. The initial system (\ref{dif-syst}), (\ref{cond-syst-1}) is the base of induction. 
The step of induction is following. 
Let us transform the system (\ref{dif-syst}), (\ref{cond-syst-1}) with constants $d_1,\ldots,d_n$ to the system with constants $\tilde{d}_1,\ldots,\tilde{d}_n$, such that $\tilde{d}_i=d_i$, $i\not=k$, $\tilde{d}_k=d_k\pm 1$, for any given $k$.\\
\\
We construct the gauge transformation as a composition of two transformations. The matrix $\hat{Y}_0=I$ from (\ref{formal-sol}) is an identity matrix, because the principle matrix $A_r$ of the system (\ref{cond-syst-1}) is diagonal. At first, let us apply the gauge transformation
$$
Y'(z)=z^{D_k^{\pm}}Y(z), 
$$   
with $D_k^{\pm}={\rm diag}(\ldots,0,\pm 1,0,\ldots)$, the $k$-th element is equal to $\pm 1$. At second, let us apply the constant transformation $Y''(z)=\hat{Y}'_0{^{-1}}Y'(z)$, where $\hat{Y}_0'$ is the first element of the formal power series (\ref{formal-sol}) of the matrix $Y'(z)$. 
It is easy to see that the matrix $A''(z)=Y''(z+1)Y''^{-1}(z)$ has the form (\ref{syst-koeff-tilde}). Actually, we can prove the lemma in a finite number of steps. As a result, we will obtain the system (\ref{syst-tilde}), (\ref{syst-koeff-tilde}) and the matrix $M(z)$ will be a composition of all the transformations at each step. (For more details of the proof, see Lemma \ref{change_sigma}.)
$\Box$\\
\\
Now we can proceed to the proof of Theorem \ref{th-BSF}.\\

{\it Proof of Theorem \ref{th-BSF}.}
Let us consider the case when Theorem \ref{th-inversprob} gives us the system (\ref{dif-syst}), (\ref{syst-koeff}) 
\begin{eqnarray}\label{syst-0}
Y(z+1)=A(z)Y(z),\qquad A(z)=A_0+\ldots+A_{r-1}z^{r-1}+A_r z^r
\end{eqnarray}
with the given monodromy $P(z)$, the given characteristic constants  $\rho_1,\ldots,\rho_n$ and shifted constants $\tilde{d}_1=d_1+l_1,\ldots,\tilde{d}_n=d_n+l_n$, $l_i\in\mathbb{Z}$, $i=1,\ldots,n$. Let us apply Lemma \ref{lemma-fixval} and construct the system (\ref{syst-tilde}), (\ref{syst-koeff-tilde}) 
\begin{eqnarray}\label{syst-1}
Y'(z+1)=A'(z)Y'(z),\qquad A'(z)=A_r' z^r+A_{r-1}'z^{r-1}+\ldots+A_0'+\ldots+A_{-s}'z^{-s}
\end{eqnarray}
with the same monodromy $P(z)$ and characteristic constants  $\rho_1,\ldots,\rho_n$ and $d_1,\ldots,d_n$. Let us denote by  $Y(z)$ and $Y'(z)$ fundamental matrices of systems (\ref{syst-0}) and (\ref{syst-1}). They have exactly the same monodromy $P(z)$. Let us denote roots of the polynomial $\det A(z)$ by $q_1,\ldots,q_{rn}$ and denote zeros and pole (z=0) of function $\det A'(z)$ by $q_1',\ldots,q_{l}'$.\\
\\
Consider the matrix $M^{-1}(z)=Y(z){Y'}^{-1}(z)$. It is easy to see that the matrix $M(z)$ is holomorphically invertible in $\mathbb{C}\setminus \{ 0\}$.
Then Sauvage's lemma (see \cite{AB}) gives us the following decomposition for matrix $M(z)$:
\begin{eqnarray}\label{sauvage}
U(z)M^{-1}(z)=z^{K}W(z),
\end{eqnarray}
where $U(z)$ is holomorphically invertible in $\mathbb{C}$, $W(z)$ is holomorphically invertible in $\overline{\mathbb{C}}\setminus\{ 0\}$, $K={\rm diag}(k_1,\ldots,k_n)$, $k_i\in\mathbb{Z}$, $i=1,\ldots,n$, $k_1\geqslant\ldots\geqslant k_n$.\\
\\
We have
$$
Y'(z)=M(z)Y(z)=W^{-1}(z)z^{-K}U(z)Y(z).
$$

Let us transform the system (\ref{syst-0}) by the transformation $\tilde{Y}(z)=U(z)Y(z)$:
\begin{eqnarray}\label{syst-01}
\tilde{Y}(z+1)=\tilde{A}(z)\tilde{Y}(z),\qquad \tilde{A}(z)=U(z+1)A(z)U^{-1}(z).
\end{eqnarray}
Let us transform the system (\ref{syst-1}) by the transformation $\tilde{Y}'(z)=W(z)Y'(z)$:
\begin{eqnarray}\label{syst-11}
\tilde{Y}'(z+1)=\tilde{A}'(z)\tilde{Y}'(z),\qquad \tilde{A}'(z)=W(z+1)A'(z)W^{-1}(z).
\end{eqnarray}
Now, it is evident that the systems (\ref{syst-01}) and (\ref{syst-11}) are connected in the following way:
$$
\tilde{Y}'(z)=z^{-K}\tilde{Y}(z).
$$
Actually, we have the following:
\begin{eqnarray}\label{calib-K}
\tilde{Y}'(z+1)\tilde{Y}'{^{-1}}(z)=\tilde{A}'(z)=(z+1)^{-K}\tilde{A}(z)z^{K}.
\end{eqnarray}
Also, the coefficient matrix $\tilde{A}(z)$ of the system (\ref{syst-01}) is a polynomial matrix of the variable $z$ (but maybe of the different degree) and the coefficient matrix $\tilde{A}'(z)$ of the system (\ref{syst-11}) is rational, but preserving the greatest power $r$ from the initial coefficient matrix $\tilde{A}(z)$ of the system (\ref{syst-1}).\\
\\
The last part of the proof will be made by induction. The pair of systems (\ref{syst-01}) and (\ref{syst-11}) is the base of induction. Now we consider the matrix $K$ from (\ref{sauvage}). Let us notice that if $L^1$-norm of $K$ is zero, the systems (\ref{syst-01}) and (\ref{syst-11}) coincide. Thus, we aim to make $\|K\|_1=\sum_{i=1}^n|k_i|=0$. To do this, we prove the following lemma.

\begin{lemma}
For the case when the difference between any two roots of $\det A(z)$ is not integer, there exist a gauge transformations of systems (\ref{syst-0}) and (\ref{syst-1}) which preserve the forms (\ref{syst-0}) and (\ref{syst-1}) of the coefficient matrices, singularities of solutions in $0$ and $\infty$ and leads to a new pair of systems with fundamental matrices 
$$
\overline{Y}'(z)=z^{-\overline{K}}\overline{Y}(z)
$$
where $\|\overline{K}\|_1\leq \|\overline{K}\|_1-1$ which means that these systems are closer than initial ones.

\end{lemma}

{\it Proof.} Let us notice that the roots of $\det A(z)$ and $\det \tilde{A}(z)$ coincide as follows from (\ref{syst-01}). 
Suppose that there exist negative numbers in $\{ k_i\}_{i\geq l}$ of $K$ so $k_1\geq\ldots\geq k_{l-1}\geq 0>k_{l}\geq\ldots\geq k_n$. Consider the points $q_1',\ldots,q_{rn}'$, such that $\det A'(q_i)=0$, $i=1,\ldots,rn$. The matrix $A'(z)$ is degenerate in every point $q_1',\ldots,q_{rn}'$. Let us consider a point $q=q_i'$ ($\det \tilde{A}(q)=0$) which means that in every point $q_i'$ there exists a nontrivial linear combination of rows of $\tilde{A}'(q_i')$ with zero sum:
\begin{eqnarray}\label{linear-comb}
\alpha_{1}\bar{a}_{j_1}(q_i')+\ldots+\alpha_{h}\bar{a}_{j_h}(q_i')=0,\qquad \alpha_{1}\cdot\ldots\cdot \alpha_{h}\not=0,
\end{eqnarray}
where $\bar{a}_j$ is a $j$-th string of the matrix $\tilde{A}'(z)$, $\alpha_j\in\mathbb{C}$, $j_1>\ldots>j_h$. 

Let us take the preliminary transformation $\hat{Y}'(z)=\left(\frac{z-q}{z}\right)^{-\overline{K}}\tilde{Y}(z)$. Now the systems 
are connected in the following way:
$$
\hat{Y}'(z)=(z-q)^{-\overline{K}}\tilde{Y}(z).
$$
Consider the case $j_h\geq l$. Take the following transformation and obtain a new pair of systems:
$$
\overline{Y}'(z)=F(z)\tilde{Y}'(z),
$$
where
\begin{eqnarray*}
F(z)=\left(
\begin{array}{ccccc}
  1 & 0 & 0 & \ldots & 0 \\
    0 & 1 & 0 & \ldots & 0\\
 \vdots &  \ddots &\ddots & \ddots & \vdots \\
   \alpha_{1}(z-q)^{k_{j_1}-k_{j_h}} & \ldots & \alpha_{h-1}(z-q)^{k_{j_{h-1}}-k_{j_h}} & \alpha_{h}  & 0 \\
  0 &  0 & \ldots & 0 & 1
\end{array}\right).
\end{eqnarray*}
This tramsformation is holomorphic in $\mathbb{C}$. Let us apply the transformation $(z-q)^{D_{j_h}^-}\overline{Y}'(z)$ and we obtain the new system pair of systems (\ref{syst-01}) and (\ref{syst-11}) with new matrix $\tilde{K}$ which is equal to $K+D_{j_h}^+$ ($\| K+D_{j_h^+}\|_1=\|k\|_{1}-1$).  

If all elements of matrix $K$ are positive then we can apply the analogous procedure on the right hand side. For more details of the proof, see Lemma \ref{distance_lemma}.

Let us prove that there is at least one linear combination (\ref{linear-comb}) in at least one of the points $q_i'$, such that $j_h\leqslant l$. Actually, let us suppose that all linear combinations (\ref{linear-comb}) in all points $q_1',\ldots,q_{rn}'$ satisfy to 
$j_h<l$. Consider the minor $\Delta(z)$ of the matrix $\tilde{A}'(z)$ formed by the last $n-l$ rows and some $n-l$ columns $\bar{b}_{i_1},\ldots,\bar{b}_{i_{n-l}}$ such that $\Delta(z)\not\equiv 0$. Such a minor exists because $\det \tilde{A}'(z)\not\equiv 0$.

In the one hand the minor $\Delta (z)$ has $rn$ roots $q_1',\ldots,q_{rn}'$. In the other hand the minor $\Delta (z)$ is a polynomial and
$$
\deg \Delta(z)\leqslant (n-l)r +k_{i_1}+\ldots+k_{i_{n-l}}-k_{l+1}-\ldots-k_n< (n-l)r.
$$
Consequently, we have a contradiction, because the polynomial $\Delta (z)\not\equiv 0$ of degree $<(n-l)r$ has $nr$ roots. (For more details of the proof, see proof of Lemma \ref{distance_lemma}.)

\medskip


$\Box$

\section{$q$-Difference Systems Case}
\subsection{Introduction to linear $q$-difference systems}
The linear $q$-difference system is the system of the form:
\begin{eqnarray}\label{syst}
Y(qz) = Q(z)Y(z), |q|>1.
\end{eqnarray}
 Here $Q(z)$ is a matrix of polynomials of degree $\mu$ or less. More general case of a rational matrix could be reduced to the polynomial by following means. Let
\begin{eqnarray}\label{rat_syst}
\hat{Y}(qz)=\hat{Q}(z)\hat{Y}(z)
\end{eqnarray} be a system with a rational matrix $\hat{Q}(z)$ and let the least common denominator of it's elements be written in the form $(z - a_1) \dots (z - a_l)$. Let also $g_i(z)$ be a solution of the following $q$-difference equation:
\begin{eqnarray}\label{one_dim_equation}
g(qz) = (z - m)g(z)
\end{eqnarray}
where $m = a_i$. Now we consider the substitution
$$
Y(z) = g_1(z) \dots g_l(z)\hat{Y}(z)
$$
in order to rewrite the system (\ref{rat_syst}) in the form (\ref{syst}):
$$
\hat{Y}(qz) = \frac{Y(qz)}{\prod_i{g_i(qz)}} = \hat{Q}(z) \frac{Y(z)}{\prod_i{g_i(z)}} \Rightarrow \frac{Y(qz)}{\prod_i{(z-a_i)g_i(z)}} = \hat{Q}(z) \frac{Y(z)}{\prod_i{g_i(z)}} \Rightarrow Y(qz) = Q(z)Y(z),
$$
$$
Q(z) = \prod_i{(z-a_i)}\hat{Q}(z)
$$
\medskip
Thus, we have obtained a polynomial coefficient matrix. It is only left to show that a solution for (\ref{one_dim_equation}) exists:
\begin{enumerate}
\item For $m = 0$ we take the transformation $t = \log_{q}{z}$ and $g(q^t) = f(t)$ and obtain:
$$
 g(qz) = zg(z) \Rightarrow f(t+1) = q^tf(t).
$$
So, we find by direct substitution that $f(t) = q^{\frac{1}{2}(t^2-t)}$ is the solution. $\Box$
\item For $m\not=0$ we take another transformation:
$$
\begin{cases}
   z = m \overline{z}
   \\
   y(z) = e^{\pi i \log_{q}{\overline{z}}} m^{\log_{q}{\overline{z}}}\overline{y}(\overline{z}),
\end{cases}
$$
Substitute it to (\ref{one_dim_equation}):
$$
e^{\pi i \log_{q}{\overline{z}}} m^{\log_{q}{\overline{z}}}\overline{y}(\overline{z})(m\overline{z} - m) = e^{\pi i (1+\log_{q}{\overline{z}})} m^{(1+\log_{q}{\overline{z}})}\overline{y}(q\overline{z}).
$$
This leads to the normal form of (\ref{one_dim_equation}):
\begin{eqnarray}\label{normal}
\overline{y}(q\overline{z}) = (1 - \overline{z})\overline{y}(\overline{z}).
\end{eqnarray}
It could (analogically to $m=0$ case) be verified by direct substitution that two solutions of (\ref{normal}) are
$$
\begin{cases}
   y_0(z) = \left(1-\frac{z}{q}\right)\left(1-\frac{z}{q^2}\right)\dots
   \\
   y_{\infty}(z) = q^{\frac{1}{2}(t^2-t)}e^{-\pi i t}\frac{1}{1 - \frac{1}{z}}\frac{1}{1-\frac{1}{qz}}\dots.
\end{cases}
$$
Both series converge since $|q|>1$.
$\Box$\\
\end{enumerate}
The fundamental existence theorems for linear $q$-difference systems (\cite{Bir-Rim}, p.561) guarantee that in general there exist two matrix solutions of (\ref{syst}):
\begin{eqnarray}\label{solution}
\begin{cases}
   Y_0(z) = A(z) z^{\begin{pmatrix}
\rho_1 & 0 & \cdots & 0 \\
0 & \rho_2 & \cdots & 0 \\
\vdots & \vdots & \ddots & \vdots \\
0 & 0 & \cdots & \rho_n
\end{pmatrix}}
   \\
   Y_{\infty}(z) = q^{\frac{\mu}{2}(t^2-t)}B(z)z^{\begin{pmatrix}
-\sigma_1 & 0 & \cdots & 0 \\
0 & -\sigma_2 & \cdots & 0 \\
\vdots & \vdots & \ddots & \vdots \\
0 & 0 & \cdots & -\sigma_n
\end{pmatrix}}
\end{cases}
\end{eqnarray}
where $A(z)$ is analytic at $z=0$, $B(z)$ is analytic at $z=\infty$ and also the determinants of the leading coefficients of series of $A(z)$ at $0$ and $B(z)$ at $\infty$ are not zero. We only consider the case when these series exist. Here \{$\rho_i$\} and \{$\sigma_j$\} are some characteristic constants which will be used later, $\mu$ is the degree of $Q(z)$, $t = \log_{q}{z}$. From (\ref{syst}), using the $q$-periodic property of the solutions we conclude that $Y_0(z)$ is analytic for $z\not=0$ and $z\not=\infty$ (because we can enlarge the vicinity of analyticity multiplying by $q$ where $|q|>1$). Similarly, $Y_{\infty}(z)$ is analytic except for $z=0$, $z=\infty$ and also poles.\\
\\
We now consider the {\it monodromy matrix} for the $q$-difference case. Here the definition is similar to the difference case, it is again the matrix which connects two solutions $Y_0$ and $Y_{\infty}$:
\begin{eqnarray}\label{def_P}
Y_0(z) = Y_{\infty}(z) P(z)
\end{eqnarray}
$P(z)$ is analytic for $z\not=0$ and $z\not=\infty$ and we can easily prove that it is $q$-periodic:
$$
\begin{cases}
Y_0(qz)=Q(z)Y_0(z)
\\
Y_{\infty}(qz)=Q(z)Y_{\infty}(z)
\\
Y_0(z) = Y_{\infty}(z) P(z)
\end{cases}
\Rightarrow Y_{\infty}^{-1}(qz)Y_0(qz) = Y_{\infty}^{-1}(z)Q^{-1}(z)Q(z)Y_0(z) \Rightarrow P(qz) = P(z). \Box
$$

\subsection{The nature of the monodromy matrix $P(z)$.}
\subsubsection{Relations for $P(z)$ elements}
We now are going to study the properties of $P(z)$. Firstly, let us take the transformation $z = q^t$ to consider our matrix $\overline{P}(t) = P(z)$ on the $t$-plane.
Let us also divide the $t$-plane into parallelograms corresponding to the periods $\omega = 1$ and $\omega' = 2\pi i/\ln{q}$. Consider $ABCD$ - one of these parallelograms with vertices $$
A = t_0, B = t_0 + 1, C = t_0 + 1 + 2\pi i/\ln{q}, D = t_0 + 2\pi i/\ln{q}.
$$
We know that $P(z) = P(qz)$, then it follows that $\overline{P}(t) = \overline{P}(t+1)$. So we have $$
\overline{P}(A) = \overline{P}(B), \overline{P}(C) = \overline{P}(D). $$\\
\\
The next step is finding the relation between $\overline{P}(A)$ and $\overline{P}(D)$, $\overline{P}(C)$ and $\overline{P}(B)$.
We recall that
\begin{eqnarray}\label{defP}
P(z) = Y_{\infty}^{-1}(z)Y_0(z).
\end{eqnarray}
If we make a positive circuit around $z = 0$, both matrices $Y_0$ and $Y_{\infty}$ will be multiplied by their monodromies (they will be branching). From the form of $Y_0$ given by (\ref{solution}) we find that $Y_0(z)$ will change to $Y_0(z)M_0$, where $$
M_0 = \begin{pmatrix}
e^{2\pi i \rho_1} & 0 & \cdots & 0 \\
0 & e^{2 \pi i\rho_2} & \cdots & 0 \\
\vdots & \vdots & \ddots & \vdots \\
0 & 0 & \cdots & e^{2\pi i\rho_n}
\end{pmatrix}.
$$
At the same time $Y_{\infty}(z)$ will change to $$
(-1)^{\mu}e^{2\pi i \mu t} e^{-2\pi^2 \mu /  \ln{q}}Y_{\infty}(z)M_{\infty}$$ where $$
M_{\infty} = \begin{pmatrix}
e^{-2\pi i\sigma_1} & 0 & \cdots & 0 \\
0 & e^{-2\pi i\sigma_2} & \cdots & 0 \\
\vdots & \vdots & \ddots & \vdots \\
0 & 0 & \cdots & e^{-2\pi i\sigma_n}
\end{pmatrix}.
$$
If we substitute both results in (\ref{defP}), we obtain that for a positive circuit around $z=0$, $P(z)$ will change to
$$
(-1)^{\mu}e^{-2\pi i \mu t} e^{2\pi^2 \mu /  \ln{q}}M_{\infty}^{-1}P(z)M_0.
$$
If $z$ makes a positive circuit around zero, $t = \log_{q}{z}$ changes to $t + 2\pi i /\ln{q}$ which corresponds to the passage across the edge $AD$ of the parallelogram $ABCD$. We finally have:
$$
\overline{P}(D) = (-1)^{\mu}e^{-2\pi i \mu t} e^{2\pi^2 \mu /  \ln{q}}M_{\infty}^{-1}\overline{P}(A)M_0,
$$
$$
\overline{P}(C) = (-1)^{\mu}e^{-2\pi i \mu t} e^{2\pi^2 \mu /  \ln{q}}M_{\infty}^{-1}\overline{P}(D)M_0.
\Box
$$

The same relations written for one element $\overline{p}_{ij}(t)$ of $\overline{P}(t)$ are the following:

\begin{eqnarray}\label{rel_p}
\begin{cases}
\overline{p}_{ij}(t + 2\pi i/\ln{q}) = (-1)^{\mu}e^{-2\pi i \mu t} e^{2\pi^2 \mu /  \ln{q}}e^{2\pi i(\sigma_i+\rho_j)} \overline{p}_{ij}(t)
\\
\overline{p}_{ij}(t+1) = \overline{p}_{ij}(t).
\end{cases}
\end{eqnarray}

\subsubsection{Explicit form of $P(z)$ elements}
We have seen that $P(z)$ is almost doubly periodic so it is natural to try to represent it using the elliptic functions. The only thing we have to be careful with is the coefficient in the first relation in (\ref{rel_p}).
\begin{proposition}[Birkhoff \cite{Bir-Rim}] 
The element $\overline{p}_{ij}(z)$ of $\overline{P}(z)$ is of the form
\begin{eqnarray}\label{exform}
\overline{p}_{ij}(t) = c_{ij}e^{\frac{-\eta \mu}{2}t^2+\left ( \eta (\sigma_i + \rho_j + v) - \eta '(\frac{ \mu }{2}+u) \right ) t}\sigma(t-a_{1}^{(i,j)})\dots \sigma(t-a_{\mu}^{(i,j)}),
\end{eqnarray}
where $t = \log_{q}{z}$, $\sigma(t)$ is the Weierstrass sigma function belonging to the periods $\omega =1$, $\omega'=\frac{2\pi i}{\ln{q}}$,
$\mu$ is the degree of polynomial $Q(z)$ in (\ref{syst}), $\{\sigma_i\}$ and $ \{\rho_j\}$ are the characteristic constants from (\ref{solution}), $u$ and $v$ are some arbitrary integers. And the following condition is satisfied:
\begin{eqnarray}\label{condition_a}
\sum_{k=1}^{\mu}{a_{k}}^{(i,j)}=\sigma_i + \rho_j - \frac{\mu\pi i}{\ln{q}} + v - \frac{2 u\pi i}{\ln{q}}. 
\end{eqnarray}
\end{proposition}

{\it Proof.}
\begin{enumerate}
\item We will search for a solution in the following form. Suppose for some element of $\overline{P}(t)$ that
\begin{eqnarray}\label{p_form}
\overline{p}_{ij}(t) = ce^{at^2 + bt}\prod_{k=1}^{\mu}{\sigma(t - a_k)}.
\end{eqnarray}
Recall the relations for Weierstrass sigma function:
\begin{eqnarray}\label{sigma}
\begin{cases}
\sigma(t+\frac{2 \pi i}{\ln{q}}) = -e^{\eta'\left (t + \frac{2 \pi i}{\ln{q}} \right )} \sigma(t)
\\
\sigma(t+1) = -e^{\eta(t+1/2)}\sigma(t)
\\
\eta\frac{2 \pi i}{\ln{q}} - \eta' = 2\pi i
\end{cases}
\end{eqnarray}
We notice that the requirement $\Re{(\omega/\omega')}>0$ is satisfied.
Let us apply the second relation from (\ref{rel_p}) to (\ref{p_form}). Using (\ref{sigma}), we easily obtain:
\begin{eqnarray}\label{ab_1}
a = -\frac{\eta\mu}{2}, b = \eta\sum_{k=1}^{\mu}a_{k} + \mu \pi i + 2 u\pi i,
\end{eqnarray}
where $u$ is an arbitrary integer. Now we take the values of $a$ and $b$ from (\ref{ab_1}) for (\ref{exform}) and substitute everything to the first relation from (\ref{rel_p}). We have then:
\begin{eqnarray}\label{ab_2}
\sum_{k=1}^{\mu}{a_{k}} + \frac{\mu \pi i}{\ln{q}} = \sigma_i + \rho_j + v - \frac{2 u\pi i}{\ln{q}},
\end{eqnarray}
where $v$ is again an arbitrary integer.
Relations (\ref{ab_1}) and (\ref{ab_2}) are equivalent to (\ref{rel_p}) written for (\ref{p_form}). Let us finally take the value for $\sum_{k=1}^{\mu}{a_k}$ from (\ref{ab_2}) and substitute it to the expression for $b$ in (\ref{ab_1}). We obtain:
$$
b = \eta(\sigma_i + \rho_j + v) - \eta'(\mu/2 + u)
$$
Afterwards we obtain:
$$
\overline{p}_{ij}(t) = c_{ij}e^{\frac{-\eta \mu}{2}t^2+\left ( \eta (\sigma_i + \rho_j + v) - \eta'(\frac{ \mu }{2}+u) \right ) t}\sigma(t-a_{1})\dots \sigma(t-a_{\mu})
$$
with the condition
$$
\sum_{k=1}^{\mu}{a_{k}} = \sigma_i + \rho_j - \frac{\mu \pi i}{\ln{q}} + v - \frac{2 u\pi i}{\ln{q}}.
$$
Here the set of $\{a_k\}$ is different for each $\overline{p}_{ij}$ so we have to write $\{a_k^{(i,j)}\}_{k=1}^\mu$.
This finishes the first part of the proof. We have obtained the characteristic constants $\{c_{ij}\}$ and $\{a_k^{(i,j)}\}_{k=1}^\mu$
which are defined up to arbitrary integers $v$ and $u$. These characteristic constants define the monodromy $P(z)$.
\item Now it is only left to prove that $\overline{p}_{ij}(t)$ could be represented in the form (\ref{p_form}). Let us take any function $\psi(t)$ which satisfies (\ref{rel_p}). Also let us say that $\phi(t)$ is the particular one of the form (\ref{p_form}). Then the function $\psi(t)/\phi(t)$ is doubly periodic analytic save for poles, and could therefore be represented as a quotient of products of sigma functions:
    $$
    \frac{\psi(t)}{\phi(t)} = C\frac{\sigma(t-\gamma_1)\dots\sigma(t-\gamma_k)}{\sigma(t-\beta_1)\dots(t-\beta_k)}, \sum\gamma_i = \sum\beta_i
    $$
If we multiply this quotient by $\phi(t)$, which is expressed as a product of sigma functions, we should obtain $\psi(t)$ - the entire function . It follows then that for each zero of $\sigma(t-\gamma_j)$ in the numerator there should exist a congruent zero $\beta_j$ of $\sigma(t-\beta_j)$ in the denominator. Such pairs of corresponding zeros may be combined leaving the coefficient $e^{ct+d}$. So, $\psi(t)$ may be also represented in the form (\ref{p_form}). This finishes the whole proof. $\Box$
\end{enumerate}

\subsection{The Riemann-Hilbert problem for $q$-difference case.}
\subsubsection{The generalized Riemann-Hilbert problem.}
Let us call the constants $\{\rho_i\}, \{\sigma_i\}, \{c_{i,j}\}, \{a_k^{(i,j)}\}_1^{\mu}$ -- the characteristic constants of the monodromy \space matrix $P(z)$ and the $q$-difference system (\ref{syst}).
 Now we are ready to formulate an analogue of the Riemann-Hilbert problem for $q$-difference systems: \\
 \\
 {\it Construct a $q$-difference system with the prescribed set of characteristic constants (i.e. with a given monodromy matrix $P(z)$ and with given power asymptotics of the solutions) and with the polynomial coefficient matrix $Q(z)$ of degree $\mu$.}\\
\\
The solution of Birkhoff is the following:
 \begin{theorem}[Birkhoff \cite{Bir-Rim}, p.566]\label{Bir_theorem}
 There exists a linear $q$-difference system (\ref{syst}) with the matrix solutions $Y_0(z), Y_\infty(z)$ either possessing prescribed characteristic constants  $\{\rho_i\}, \{\sigma_i\}, \{c_{i,j}\}, \{a_k^{(i,j)}\}_1^{\mu}$ or else constants  $\{\rho_j\}, \{\sigma_j+l_j\}, \{c_{i,j}\}, \{a_k^{(i,j)}\}_1^{\mu}$, where $\{l_j\}$ are integers. For an arbitrary loop about $z=0$ which cuts each spiral
 \begin{eqnarray}\label{spiral}
 \theta = c + \frac{\arg(q)}{\ln{|q|}}\ln{r},\text{   $r$,$\theta$ - polar coordinates}
 \end{eqnarray}

 only once and does not pass through points $z$ such that $\det P(z)=0$, there exist matrices $Y_0(z), Y_\infty(z)$ with the further property that $\det Y_0(z)\not=0$ within or along the loop while the elements of $Y_\infty(z)$ are analytic and $\det Y_\infty(z)$ is not zero without the loop.
 \end{theorem}
\subsubsection{Refined version of Theorem \ref{Bir_theorem}}
The solution of Birkhoff could be insufficient for the case when it derives the system with shifted constants $\{\sigma_j+l_j\}$. Let us develop further this result and improve Theorem \ref{Bir_theorem} in the following way:
\begin{theorem}\label{main}
There exists a linear $q$-difference system (\ref{syst}) with the matrix solutions $Y_0(z), Y_\infty(z)$ possessing the prescribed set of characteristic constants  $\{\rho_i\}, \{\sigma_i\}$ and $\{c_{i,j}\}, \{a_k^{(i,j)}\}_1^{\mu}$ such that the monodromy matrix $P(z)$ does not have multiple zeros. For an arbitrary loop about $z=0$ which cuts each spiral (\ref{spiral}) only once and does not pass through the points $z$ such that $\det P(z)=0$, there exist matrices $Y_0(z), Y_\infty(z)$ with the further property that $\det Y_0(z)\not=0$ within or along the loop while the elements of $Y_\infty(z)$ are analytic and $\det Y_\infty(z)$ is not zero without the loop.

\end{theorem}

We will prove this theorem in several steps. Firstly, we have to prove the following lemma.
\begin{lemma}\label{change_sigma}
Let us consider a linear $q$-difference system (\ref{syst}) with the polynomial matrix $Q(z)=Q_{\mu}z^\mu+\dots+Q_0$ of degree $\mu$ and with the characteristic constants $\{\rho_j\}, \{\sigma_j+l_j\}, \{c_{i,j}\}, \{a_k^{(i,j)}\}_1^{\mu}$, where $\{l_j\}$ are integers. Then, for the case when the difference between any two nulls of $\det Q(z)$ is not integer, there exists a linear $q$-difference system $Y'(qz)=Q'(z)Y'(z)$ with the characteristic constants  $\{\rho_i\}, \{\sigma_i\}, \{c_{i,j}\}, \{a_k^{(i,j)}\}_1^{\mu}$ and with a rational matrix $Q'(z)=Q'_\mu z^\mu + \ldots + Q'_0 + \ldots Q_{-s}z^{-s}$ which satisfies the following relation:
$$
Q'(z) = M(qz)Q(z)M^{-1}(z)
$$
for some rational matrix $M(z)$.
\end{lemma}
{\it Proof.}
Firstly, let us notice that we can consider the leading coefficient of $Y_\infty$ as $B_0 = I$ -- the identity matrix since $B_0$ is constant and $\det B_0\not= 0$. Also, we can consider the coefficient $Q_\mu$ in $Q(z)$ in the system (\ref{syst}) as $Q_\mu = q^{{\rm diag}(-\sigma_1,\ldots,-\sigma_n)}$. Let us prove these facts.\\
\\
Let us have an arbitrary system $\hat{Y}(qz)=\hat{Q}(z)\hat{Y}(z)$ which satisfies the conditions of the lemma. We will now show that we can easily transform it to the system $Y(qz)=Q(z)Y(z)$ with $Q_\mu =  q^{{\rm diag}(-\sigma_1,\ldots,-\sigma_n)}$ and $B_0=I$. We take the transformation $Y(z) = \hat{B}_0^{-1} \hat{Y}(z)$ and obtain at once that $B_0=I$ in $Y_\infty$. Now let us consider $Q(z)$ and use the explicit form (\ref{solution}) of $Y_\infty$:
 $$
 Q(z) = Y_{\infty}(qz)Y_{\infty}^{-1}(z)=q^{\mu t}B(qz)(qz)^{{\rm diag}(-\sigma_1,\ldots,-\sigma_n)}B^{-1}(z)z^{{\rm diag}(\sigma_1,\ldots,\sigma_n)}=
 $$$$=z^{\mu}\left(B_0+\frac{B_1}{qz}+\dots\right)(qz)^{{\rm diag}(-\sigma_1,\ldots,-\sigma_n)}{\left(B_0+\frac{B_1}{z}+\dots\right)}^{-1}z^{{\rm diag}(\sigma_1,\ldots,\sigma_n)}.
$$ 
Now we want to find the coefficient of the greatest power $\mu$ in polynomial $Q(z)$ so we have to consider only the leading coefficients in the series. Using the fact that the leading coefficient of the $B^{-1}(z)$ is $B_0^{-1}=B_0=I$ we obtain:
$$
Q_{\mu} = B_0(qz)^{{\rm diag}(-\sigma_1,\ldots,-\sigma_n)}B_0^{-1}z^{{\rm diag}(\sigma_1,\ldots,\sigma_n)}=q^{{\rm diag}(-\sigma_1,\ldots,-\sigma_n)}.
$$

Now we can get back to the proof of the whole Lemma \ref{change_sigma}.\\
\\
Let us consider the transformations which changes one particular $l_k$ of the $\{l_j\}$ to $l_k\pm1$. This will finish the proof in a finite number of steps.\\
\\
We take the following transformations:
\begin{enumerate}
\item $Y''(z) = z^{D_k^{\pm}} Y(z)$
where $D_k^{\pm} = {\rm diag}(\ldots,0,\pm 1,0,\ldots)$ with the $k$-th element equal to $\pm 1$
\item $\tilde{Y}(z) = B_0''^{-1}Y''(z),$
\end{enumerate}
Let us study how transformation 1 will affect the solutions $Y_0$ and $Y_{\infty}$. At first, we consider $Y_{\infty}$:
$$
z^{D_k^{\pm}}Y_{\infty} = q^{\frac{\mu}{2}(t^2-t)}z^{D_k^{\pm}}\left(B_0 + \frac{B_1}{z}+\dots\right)z^{{\rm diag}(\sigma_1+l_1, \dots, \sigma_k+l_k,\dots, \sigma_n+l_n)} =
$$$$ = q^{\frac{\mu}{2}(t^2-t)}\left(B_0'' + \frac{B_1''}{z}+\dots\right)z^{D_k^{\pm}}z^{{\rm diag}(\sigma_1+l_1, \dots, \sigma_k+l_k,\dots, \sigma_n+l_n)}=
$$$$ = q^{\frac{\mu}{2}(t^2-t)}\left(B_0'' + \frac{B_1''}{z}+\dots\right)z^{{\rm diag}(\sigma_1+l_1, \dots, \sigma_k+l_k\pm1,\dots, \sigma_n+l_n)}
$$
So we have to carry our $z^{D_k^{\pm}}$ through the series $B(z)$ in order to obtain $z^{D_k^{\pm}}$ on the right of $B(z)$. We also have to assure that a new series will not have negative powers of $\frac{1}{z}$. It is evident that the only coefficient we have to be careful about is $B_0=I$. It is obvious that $z^{D_k^{\pm}}B_0$ just multiplies $k$-th row of $B_0$ by $z^{\pm1}$, while $B_0z^{D_k^{\pm}}$ multiplies $k$-th column of $B_0$ by $z^{\pm1}$. Let us consider two cases -- of $z^{D_k^{+}}$ and of $z^{D_k^{-}}$ and see that in both of them the series will remain analytic.
\begin{enumerate}
\item The first case: $z^{D_k^{-}}$. Let us call $z^{D_k^{-}}B(z) = \hat{B}(z)$.
$$z^{D_k^{-}}B(z) = z^{D_k^{-}}B_0+z^{D_k^{-}}\frac{B_1}{z}+\dots= \bordermatrix{
& & &  &\bold{k} \cr
& 1 & 0 &\dots & 0 & \dots & 0  \cr
& 0 & 1 & \dots & \vdots & \dots & 0\cr
 & \vdots  &  & \ddots & \vdots& &\vdots \cr
\bold{k} & 0 & \dots & \dots & \frac{1}{z} & \dots &0\cr
& \vdots &  &  & \vdots &\ddots & \vdots \cr
& 0 & \dots & \dots & 0 &\dots & 1\cr}+z^{D_k^{-}}\frac{B_1}{z}+\dots=
$$$$
=\bordermatrix{
& & &  &\bold{k} \cr
& 1 & 0 &\dots & 0 & \dots & 0  \cr
& 0 & 1 & \dots & \vdots & \dots & 0\cr
 & \vdots  &  & \ddots & \vdots& &\vdots \cr
\bold{k} & 0 & \dots & \dots & 0 & \dots &0\cr
& \vdots &  &  & \vdots &\ddots & \vdots \cr
& 0 & \dots & \dots & 0 &\dots & 1\cr}z^{D_k^{+}}z^{D_k^{-}}+\frac{\hat{B}_1}{z}z^{D_k^{+}}z^{D_k^{-}}+\dots=B''_0z^{D_k^{-}}+\frac{B_1''}{z}z^{D_k^{-}}+\dots=
$$$$
=\left(\bordermatrix{
& & &  &\bold{k} \cr
& 1 & 0 &\dots & * & \dots & 0  \cr
& 0 & 1 & \dots & * & \dots & 0\cr
 & \vdots  &  & \ddots & *& &\vdots \cr
\bold{k} & 0 & \dots & \dots & 1 & \dots &0\cr
& \vdots &  &  & * &\ddots & \vdots \cr
& 0 & \dots & \dots & * &\dots & 1\cr}+\frac{B_1''}{z}+\dots\right)z^{D_k^{-}}
$$
\item The second case: $z^{D_k^{+}}$. Let us call $z^{D_k^{+}}B(z) = \hat{B}(z)$.
$$z^{D_k^{+}}B(z) = z^{D_k^{+}}B_0+z^{D_k^{+}}\frac{B_1}{z}+\dots= \bordermatrix{
& & &  &\bold{k} \cr
& 1 & 0 &\dots & 0 & \dots & 0  \cr
& 0 & 1 & \dots & \vdots & \dots & 0\cr
 & \vdots  &  & \ddots & \vdots& &\vdots \cr
\bold{k} & 0 & \dots & \dots & z & \dots &0\cr
& \vdots &  &  & \vdots &\ddots & \vdots \cr
& 0 & \dots & \dots & 0 &\dots & 1\cr}+z^{D_k^{-}}\frac{B_1}{z}+\dots=
$$$$
=\bordermatrix{
&  &  &\bold{k} \cr
& 0  & \dots & 0 & \dots & 0\cr
 & \vdots    & \ddots & \vdots& &\vdots \cr
\bold{k} & 0  & \dots & z & \dots &0\cr
& \vdots &    & \vdots &\ddots & \vdots \cr
& 0  & \dots & 0 &\dots & 0\cr}z^{D_k^{-}}z^{D_k^{+}}+\bordermatrix{
& &   &\bold{k} \cr
& 1 & \dots & 0 & \dots & 0\cr
 & \vdots    & \ddots & \vdots& &\vdots \cr
\bold{k} & * & * & * & * &*\cr
& \vdots   &  & \vdots &\ddots & \vdots \cr
& 0  & \dots & 0 &\dots & 1\cr}z^{D_k^{-}}z^{D_k^{+}}+\frac{\hat{B}_1}{z}z^{D_k^{-}}z^{D_k^{+}}+\dots=
$$$$
=0 + B''_0z^{D_k^{+}}+\frac{B_1''}{z}z^{D_k^{+}}+\dots=\left(\bordermatrix{
& & &  &\bold{k} \cr
& 1 & 0 &\dots & 0 & \dots & 0  \cr
& 0 & 1 & \dots & \vdots & \dots & 0\cr
 & \vdots  &  & \ddots & \vdots& &\vdots \cr
\bold{k} & * & * & * & 1 & * &*\cr
& \vdots &  &  & \vdots &\ddots & \vdots \cr
& 0 & \dots & \dots & 0 &\dots & 1\cr}+\frac{B_1''}{z}+\dots\right)z^{D_k^{+}}
$$
\end{enumerate}
Afterwards we obtain a new series $B''(z)$ in $Y_\infty''(z)$ and $z^{D_k^{\pm}}$ on the right of the series. Then we apply $z^{D_k^{\pm}}$ to $z^{{\rm diag}(\sigma_1+l_1, \dots, \sigma_k+l_k,\dots, \sigma_n+l_n)}$ and change the $l_k$ to $l_k\pm1$. Let us notice that $B''_0\not=I$ and $\det B''_0 =1$ (it is evident from the explicit form of $B''_0$). We apply the second transformation in order to make the leading coefficient an identity matrix again.
Let us notice that these transformations will probably change the properties of $Y_0(z)$ and could increase the number of nulls or poles but only in $z=0$ and $z=\infty$ because the first transformation has singularities in these points and the second is constant. But the order of nulls and poles of $Y''_\infty$ will remain correct.\\
\\
One step leads us to the system $\tilde{Y}(qz) = \tilde{Q}(z)\tilde{Y}(z)$ for which we have $$
\tilde{Q}(z) = B_0''^{-1}(qz)^{D_k^{\pm}}Q(z)(z^{D_k^{\pm}})^{-1}B_0''.
$$
Let us prove that the greatest power in rational matrix $\tilde{Q}(z)$ is $\mu$ so $\tilde{Q}(z) = \tilde{Q}_\mu z^\mu + \dots + \tilde{Q}_0+\dots+\tilde{Q}_{-1}z^{-1}$. It is obvious that $B_0''$ does not affect the power since $B_0''$ is constant. Thus, it is sufficient to prove that $(qz)^{D_k^{\pm}}Q(z)(z^{D_k^{\pm}})^{-1}$ preserves the greatest power $\mu$. We again have to consider two cases and use the properties of $z^{D_k^{\pm}}$.
\begin{enumerate}
\item $z^{D_k^{+}}$:
$$
(qz)^{D_k^{+}}Q(z)(z^{D_k^{+}})^{-1}=(qz)^{D_k^{+}}(Q_{\mu}z^\mu+\dots+Q_0)(z^{D_k^{+}})^{-1}=
$$$$
=(Q'''_{\mu+1}z^{\mu+1}+Q'''_{\mu}z^{\mu}+\dots+Q'''_0)z^{D_k^{-}}=Q''_\mu z^\mu+\dots+Q''_0+Q''_{-1}z^{-1}
$$since$$
Q_\mu = {\rm diag}(*\dots*) \Rightarrow Q'''_{\mu+1} = {\rm diag}(0\ldots,0,*,0,\ldots0).
$$
\item $z^{D_k^{-}}$:
$$
(qz)^{D_k^{-}}Q(z)(z^{D_k^{-}})^{-1}=(qz)^{D_k^{-}}(Q_{\mu}z^\mu+\dots+Q_0)z^{D_k^{+}}=
$$$$
=(qz)^{D_k^{-}}(Q'''_{\mu+1}z^{\mu+1}+Q'''_{\mu}z^{\mu}+\dots+Q'''_0)=Q''_\mu z^\mu+\dots+Q''_0+Q''_{-1}z^{-1}
$$since$$
Q_\mu = {\rm diag}(*\dots*) \Rightarrow Q'''_{\mu+1} = {\rm diag}(0\ldots,0,*,0,\ldots0).
$$
\end{enumerate}
So we have changed $l_k$ to $l_k\pm1$. Thus, we can finish the entire proof in a finite number of steps.
 As we have already seen, one step -- the composition of the two latter transformations is rational as a composition of rational and constant matrices which preserve the greatest power of the coefficient matrix. Thus, the finite number of steps will lead to the rational transformation which preserves the greatest power too. Let us call this composition $M(z)$ which will have singularities only in $0$ and $\infty$. So, we finally make the transformation $Y'(z) = M(z)Y(z)$ and obtain the system $Y'(qz)=Q'(z)Y'(z)$ with the characteristic constants  $\{\rho_i\}, \{\sigma_i\}, \{c_{i,j}\}, \{a_k^{(i,j)}\}_1^{\mu}$ and with the rational matrix $Q'(z)$ which satisfies the following relation:
$$
Q'(z) = M(qz)Q(z)M^{-1}(z). \; \Box
$$

{\it Proof of the Theorem \ref{main}}.

 Let us have the prescribed set of characteristic constants $\{\rho_i\}$, $\{\sigma_i\}$, $\{c_{i,j}\}$, $\{a_k^{(i,j)}\}_1^{\mu}$. Let us apply Theorem \ref{Bir_theorem}. And let us construct a system
    \begin{eqnarray}\label{system} 
    Y(qz) = Q(z)Y(z)
    \end{eqnarray}
    from Theorem \ref{Bir_theorem}. The case when this system has precisely the prescribed set of constants is not interesting since there is nothing to prove here. Let us consider the case when the system from Theorem \ref{Bir_theorem} possesses the characteristic constants $\{\rho_j\}, \{\sigma_j+l_j\}, \{c_{i,j}\}, \{a_k^{(i,j)}\}_1^{\mu}$. The matrix $Q(z)$ is of the form
    $$
    Q(z) = Q_\mu z^\mu+\dots+Q_0. 
    $$
    Let us denote $\alpha_1\dots\alpha_{\mu n}$ the roots of the polynomial $\det Q(z)$.\\
\\
We apply Lemma \ref{change_sigma} and obtain a system 
\begin{eqnarray} \label{system'}
Y'(qz) = Q'(z)Y(z)
\end{eqnarray}
where
$$
Q'(z) = M(qz)Q(z)M^{-1}(z),  \text{\space}  Q'(z)=Q'_\mu z^\mu + \dots + Q'_0 + \dots Q_{-s}z^{-s}.
$$
As we know, $M(z)$ is rational and holomorphically invertible in $\mathbb{C}\backslash \{0\}$. Thus, we can apply Sauvage's Lemma and obtain the following decomposition:
$$
U(z)M^{-1}(z)=z^{-D}W(z), 
$$
where $U(z)$ is holomorphically invertible in $\mathbb{C}$, $W(z)$ is holomorphically invertible in $\overline{\mathbb{C}}\backslash \{0\}$ and $D={\rm diag}(d_1\dots d_n)$ such that $d_1\ge d_2 \ge \dots \ge d_n$ where $\{d_i\}$ are integers.\\
\\
Let us transform both systems -- the initial (\ref{system}) and the improved (\ref{system'}) at the same time. Let us emphasize that the system (\ref{system'}) is "good" in infinity (because it has the correct characteristic constants in $Y_\infty$) and the system (\ref{system}) is "good" in zero (since it has the correct properties of $Y_0$). To make our notations more clear, we notice that the system (\ref{system'}) and all the the systems with prime that we are now going to construct correspond to the "good behavior" in infinity, and, vice versa, the system (\ref{system}) and all the other systems without prime correspond to zero. Now, consider a new pair of systems:
\begin{eqnarray} \label{init_tr}
\tilde{Y}(qz) = \tilde{Q}(z)\tilde{Y}(z), \text{\space  } \tilde{Q}(z) = U(qz) Q(z) U^{-1}(z) 
\end{eqnarray}
where $\tilde{Y}(z)=U(z)Y(z)$ and 
\begin{eqnarray} \label{new_tr}
\tilde{Y}'(qz) = \tilde{Q}'(z)\tilde{Y}'(z), \text{\space  } \tilde{Q}'(z) = W(qz) Q'(z) W^{-1}(z) 
\end{eqnarray}
where $\tilde{Y}'(z)=W(z)Y'(z)$.
Matrices $U(z)$ and $W(z)$ are holomorphically invertible in $\mathbb{C}$ and $\overline{\mathbb{C}}\backslash\{0\}$ respectively, thus 
\begin{eqnarray} \label{form_Q}
 \tilde{Q}(z) = \tilde{Q}_{\mu_1}z^{\mu_1} + \dots + \tilde{Q}_0, \text{   } \tilde{Q}'(z) = \tilde{Q}'_{\mu}z^{\mu} + \dots + \tilde{Q}'_{-s_1}z^{-s_1}
\end{eqnarray}
for some $\mu_1$ and $s_1$.
The systems (\ref{init_tr}) and (\ref{new_tr}) are connected with the gauge transformation:
\begin{eqnarray} \label{gauge}
\tilde{Y}'(z)=z^{D}\tilde{Y}(z)
\end{eqnarray}
which leads to 
\begin{eqnarray}\label{q_gauge}
\tilde{Q}'(z) = (qz)^{D} \tilde{Q}(z) z^{-D}.
\end{eqnarray}
From (\ref{q_gauge}) and the form (\ref{form_Q}) of matrices $\tilde{Q}'(z)$ and $\tilde{Q}(z)$ we conclude that the elements $\tilde{q}_{ij}(z)$ of $\tilde{Q}(z)$ are polynomials of degree $\deg{\tilde{q}_{ij}(z)} \le \mu+d_j-d_i$.\\
\\
The following is the idea of the transformations as an illustration of the process. We simultaneously transform the systems corresponding to zero and to infinity and afterwards we make them almost coincident and connected through (\ref{gauge}). 
$$0:
\begin{cases}
  Y(qz) = Q(z)Y(z) 
\\
 Q(z)=Q_{\mu}z^\mu+\dots+Q_0 
\end{cases}
\stackrel{U}{\longrightarrow}
\begin{cases}
\tilde{Y}(qz) = \tilde{Q}(z)\tilde{Y}(z)
\\
\tilde{Q}(z) = \tilde{Q}_{\mu_1}z^{\mu_1} + \dots + \tilde{Q}_0
\end{cases}
$$
$$\infty:
\begin{cases}
  Y'(qz) = Q'(z)Y'(z)
\\ 
 Q'(z)=Q'_\mu z^\mu + \dots + Q'_{-s}z^{-s}
\end{cases}
\stackrel{W}{\longrightarrow}
\begin{cases}
\tilde{Y}'(qz) = \tilde{Q}'(z)\tilde{Y}'(z)
\\
\tilde{Q}'(z) = \tilde{Q}'_{\mu}z^{\mu} + \dots + \tilde{Q}'_{-s_1}z^{-s_1}
\end{cases},
$$
Now let us notice that if the $L^1$-norm of the matrix $D$ from (\ref{gauge}) equals to zero, the systems coincide. So we have to make $\|D\|_1$ zero in order to finish the proof of Theorem \ref{main} or, which is equivalent, we have to transform the matrix $z^D$ to an identity matrix. To guarantee that it is possible, we will prove the following lemma.

\begin{lemma} \label{distance_lemma}
For the case when the difference between any two nulls of $\det Q(z)$ is not integer, there exist transformations of systems (\ref{init_tr}) and (\ref{new_tr}) which preserve the form (\ref{form_Q}) of the coefficient matrices, singularities of solutions in $0$ and $\infty$ and lead to a new pair of systems with fundamental matrices $\overline{Y}(z)$ and $\overline{Y}'(z)$ which are connected with the following gauge transformation:
$$
\overline{Y}'(z) = z^{\overline{D}}\overline{Y}(z)
$$
where $\|\overline{D}\|_1 \le \|D\|_1 - 1$ which means that these systems are closer than the initial ones.
\end{lemma}
{\it Proof.} 
\begin{enumerate}
\item We firstly notice that the roots of $\det{Q(z)}$ and $\det{\tilde{Q}(z)}$ coincide as follows from (\ref{init_tr}). Let us consider one point $\alpha=\alpha_i$ for which $\det{\tilde{Q}(z)}=0$.\\
\\
Suppose that there exist negative numbers in $\{d_i\}$ of $D$ so $d_1\ge \ldots \ge 0 >  d_k \ge \ldots \ge d_n.$
Now $\det{\tilde{Q}(\alpha)}=0$ which means that there exists a nontrivial linear combination with zero sum of rows $\tilde{q}_{i_j}$ of $\tilde{Q}(z)$ where $j=1,\dots,l$ with non-zero $\{s_i\}$ for some rows indexes $i_1 \le\dots\le i_l$:
\begin{eqnarray}\label{combo}
s_1\tilde{q}_{i_1}(\alpha)+\ldots+s_l\tilde{q}_{i_l}(\alpha)=0.
\end{eqnarray}
Let us take the preliminary transformation $\hat{Y}'(z) = \left(\frac{z-\alpha}{z}\right)^{D}\tilde{Y}'(z).$ Now the systems are connected in the following way:$$
\hat{Y}'(z)=(z-\alpha)^{D}\tilde{Y}(z).
$$
Consider the case when $i_l \ge k$ (i.e. there exists a negative number in $D$ corresponding to the non-zero $s_l$). We take the following transformation and obtain a new pair of systems:
$$
\overline{Y}(z) = (z-\alpha)^{D^{-}_{i_l}}H\tilde{Y}(z), \text{   } \overline{Y}'(z) = (z-\alpha)^{D}H(z-\alpha)^{-D}\hat{Y}'(z)
$$
where $H$ is a constant matrix of the form $H=I+H'$. All elements of $H'$ are zeros except for $h'_{i_li_j}=s_j/s_l, j=1,\ldots,l-1$. The matrix $D^{-}_{i_l}$ is defined in Lemma \ref{change_sigma}. The matrix $(z-\alpha)^DH(z-\alpha)^{-D}$ is a polynomial of the variable $\frac{1}{z-\alpha}$ since $H$ is lower triangular by definition and the diagonal elements of $D$ decrease. The coefficient matrix will transform in the following way:$$
\overline{Q}'(z) = (qz-\alpha)^{D}H(qz-\alpha)^{-D}\hat{Q}'(z)(z-\alpha)^{D}H^{-1}(z-\alpha)^{-D}.
$$Thus, the greatest power of $\overline{Q}'$ would not increase and would remain $\mu$ (since $(z-\alpha)^{D}H(z-\alpha)^{-D}$ is holomorphically invertible in infinity). It is also evident that:
$$
\overline{Q}(z) =  (qz-\alpha)^{D^{-}_{i_l}}H\tilde{Q}(z)H^{-1}(z-\alpha)^{D^{+}_{i_l}}
$$
is holomorphic in $\mathbb{C}$ because $i_l$-th row of $H\tilde{Q}(z)$ is zero and division by $(qz-\alpha)$ would not create a singularity. Now for matrices $\overline{Y}(z)$ and $\overline{Y}'(z)$ we can write:
$$
\overline{Y}'(z)=(z-\alpha)^{D-D^{-}_{i_l}}\overline{Y}(z), \text{   } \overline{D} = D-D^{-}_{i_l}.
$$
Now we have to apply the inverse of the preliminary transformation to make the systems be connected in the correct way (we have to emphasize that the preliminary transformation and its inverse is holomorphically invertible in infinity and is applied only to the system corresponding to infinity):
$$
\hat{\overline{Y}}'(z) = \left(\frac{z}{z-\alpha}\right)^{\overline{D}} \overline{Y}'(z).
$$
So, we have added a unit to one of the negative elements of $D$. Thus, $\|\overline{D}\|_1 \le \|D\|_1 - 1$. It is only left to return the descending order to the elements of $\overline{D}$ which could be easily achieved applying a permutation matrix which is a constant transformation. \\
\\
\item Consider the case when the linear combination (\ref{combo}) satisfy $i_l<k$. Let us show that the case when for all $\alpha=\alpha_j, j = 1,\dots,\mu n$ the linear combination (\ref{combo}) satisfies $i_l<k$ is not possible. Let us prove by contrapositive. We suppose that all $\{\alpha_i\}$ are different and do not differ by an integer. Then in every point $\alpha_i$ initial $k-1$ rows of $\tilde{Q}(\alpha_i)$ are linearly dependent. Then, let us take a nontrivial minor $\Delta(z)$ which is formed with initial $k-1$ rows and some $k-1$ columns and which obviously exists. Then, using the (\ref{q_gauge}) and the fact that the elements of $D$ decrease, we obtain that the degree of such a minor satisfies the following:
$$ 
\deg(\Delta(z)) \le \mu(k-1) - (d_1 + \ldots + d_k-1) + (d_{i_1} + \ldots + d_{i_{k-1}}) \le \mu(k-1).
$$
On the other hand, the roots of $\Delta(z)$ include all $\alpha_1, \dots, \alpha_{\mu n}$. Thus, we have a contradiction since the number of the roots of a polynomial could not exceed its degree.\\
\\
\item Now, we consider the case when there are no negative $d_i$. It could be proven analogically, but we have to act on columns, not on rows. Using the previous part of the proof, we could eliminate all the negative numbers in $\{d_i\}$ of $D$. If there are no positive numbers, then the proof is finished, otherwise there should be some positive elements of $D$, so $d_1\ge \ldots \ge d_k > 0 \ge \ldots \ge d_n.$
Now, $\det{\tilde{Q}(\alpha)}=0$ which means that there exists a nontrivial linear combination with zero sum of columns $\tilde{q}_{i_j}$ of $\tilde{Q}(z)$ where $j=1,\dots,l$ with non-zero $\{s_i\}$ for some columns indexes $i_1 \le\dots\le i_l$:
\begin{eqnarray}\label{combo}
s_1\tilde{q}_{i_1}(\alpha)+\ldots+s_l\tilde{q}_{i_l}(\alpha)=0.
\end{eqnarray}
Let us take the same preliminary transformation $\hat{Y}'(z) = \left(\frac{z-\alpha}{z}\right)^{D}\tilde{Y}'(z).$ Then:$$
\hat{Y}'(z)=(z-\alpha)^{D}\tilde{Y}(z).
$$
Consider the case when $i_1 \le k$ (i.e. there exists a positive number in $D$ corresponding to the non-zero $s_1$). Now, we take the following transformations and obtain a new pair of systems:
$$
\overline{Y}(z) = (z-\alpha)^{D^{+}_{i_1}}H^{-1}\tilde{Y}(z), \text{   } \overline{Y}'(z) = (z-\alpha)^{D}H^{-1}(z-\alpha)^{-D}\hat{Y}'(z)
$$
where $H$ is a constant matrix of the form $H=I+H'$. All elements of $H'$ are zeros except for $h'_{i_j i_1}=s_j/s_1, j=1,\ldots,l-1$ (in this case, we put the coefficients of the linear combination (\ref{combo}) to the $i_1$-column of $H'$). We will now prove that the matrix $(z-\alpha)^DH^{-1}(z-\alpha)^{-D}$ is a polynomial of the variable $\frac{1}{z-\alpha}$. We have: $(z-\alpha)^DH^{-1}(z-\alpha)^{-D} = ((z-\alpha)^DH(z-\alpha)^{-D})^{-1}$ where $(z-\alpha)^DH(z-\alpha)^{-D}$ is a polynomial of $\frac{1}{z-\alpha}$ as was proven before. Since $\det (z-\alpha)^DH(z-\alpha)^{-D} = 1$, we have by Cramer's rule that the inverse matrix $(z-\alpha)^DH^{-1}(z-\alpha)^{-D}$ is a polynomial of $\frac{1}{z-\alpha}$ too. Also, the coefficient matrix will transform in the following way:$$
\overline{Q}'(z) = (qz-\alpha)^{D}H^{-1}(qz-\alpha)^{-D}\hat{Q}'(z)(z-\alpha)^{D}H(z-\alpha)^{-D}.
$$Thus, the greatest power of $\overline{Q}'$ would not increase and would remain $\mu$ (since $(z-\alpha)^{D}H^{-1}(z-\alpha)^{-D}$ is holomorphically invertible in infinity). It is also evident that:
$$
\overline{Q}(z) =  (qz-\alpha)^{D^{+}_{i_1}}H^{-1}\tilde{Q}(z)H(z-\alpha)^{D^{-}_{i_1}}
$$
is holomorphic in $\mathbb{C}$ because $i_1$-th column of $\tilde{Q}(z)H$ is zero and division by $(z-\alpha)$ would not create a singularity. Now for the matrices $\overline{Y}(z)$ and $\overline{Y}'(z)$ we can write:
$$
\overline{Y}'(z)=(z-\alpha)^{D-D^{+}_{i_1}}\overline{Y}(z), \text{   } \overline{D} = D-D^{+}_{i_1}.
$$
Now we have to apply the inverse of the preliminary transformation to make the systems be connected in the correct way:
$$
\hat{\overline{Y}}'(z) = \left(\frac{z}{z-\alpha}\right)^{\overline{D}} \overline{Y}'(z).
$$
So, we have subtracted a unit from one of the positive elements of $D$. Thus, $\|\overline{D}\|_1 \le \|D\|_1 - 1$. Again, it is left to return the descending order to the elements of $\overline{D}$ which could be easily achieved applying a permutation matrix which is a constant transformation. Also, it is needed to assure that there exists a positive (non-zero) element of $D$ corresponding to $i_1$ and the proof is literally the same as the proof of the similar fact in the second point. $ \Box$ \\
\end{enumerate}
Let us return to the proof of Theorem \ref{main}. Applying Lemma \ref{distance_lemma} finitely many times while $\|D\|_1>0$ we transform our systems (\ref{init_tr}) and (\ref{new_tr}) to a new pair of systems:
\begin{eqnarray} \label{final_zero}
\hat{\overline{Y}}(qz)=\hat{\overline{Q}}(z)\hat{\overline{Y}}(z), \text{   } \hat{\overline{Q}}(z) = \hat{\overline{Q}}_{\mu_1}z^{\mu_1} + \dots + \hat{\overline{Q}}_0, 
\end{eqnarray}
\begin{eqnarray} \label{final_infty}
\hat{\overline{Y}}'(qz)=\hat{\overline{Q}}'(z)\hat{\overline{Y}}'(z), \text{   } \hat{\overline{Q}}'(z) = \hat{\overline{Q}}'_{\mu}z^{\mu} + \dots + \hat{\overline{Q}}'_0 + \dots + \hat{\overline{Q}}'_{-s_1}z^{-s_1}
\end{eqnarray}
for some $\mu_1$ and $s_1$.
We made $\|D\|_1=0$ which means that systems (\ref{final_zero}) and (\ref{final_infty}) coincide so we have one new system possessing all the necessary properties and the prescribed set of characteristic constants:
\begin{eqnarray}\label{good_final}
\hat{\overline{Y}}'(qz)=\hat{\overline{Q}}'(z)\hat{\overline{Y}}'(z), \text{   } \hat{\overline{Q}}'(z) = \hat{\overline{Q}}'_{\mu}z^{\mu} + \dots + \hat{\overline{Q}}'_0.
\end{eqnarray}
The system (\ref{final_infty}) had the correct constants $\{\sigma_i\}$, the system ($\ref{final_zero}$) had the correct set $\{\rho_j\}$. Also, all the transformations we made preserve the monodromy matrix $P(z)$ and thus the sets of constants responsible for the monodromy is also correct in our final system. More importantly, the monodromy itself of the final system (\ref{good_final}) just coincides with that of the initial system (\ref{system}):$$
\hat{\overline{P}}'(z)=P(z).$$ In addition, we emphasize that the condition (\ref{condition_a}) holds since it is defined up to arbitrary integers $u$ and $v$ while we made the integer shifts of $\{\sigma_i\}$. Thus, we have constructed a $q$-difference system (\ref{good_final}) with the prescribed set of characteristic constants, the given monodromy and the coefficient matrix of degree $\mu$.
$\Box$

\section{Acknowledgments}

The authors are grateful to Renat Gontsov for useful discussions. The work of the first author was supported by grants RFBR 17-01-00515-a, RFBR-CNRS 16-51-150005-cnrs-a

Vyugin I.V.\\
Insitute for Information Transmission Problems RAS,  \\
and\\
National Research University Higher School of Economics,\\
{\it vyugin@gmail.com}.\\
\\
Levin R.I.\\
National Research University Higher School of Economics,\\
{\it romalevi@yandex.ru}.\\

\end{document}